\documentclass[final,3p,times]{elsarticle}

\usepackage{graphicx}

\usepackage{amssymb}

\biboptions{comma,square}

\def\epsilon{\varepsilon}
\def\bigoh{\mathcal{O}}
\def\0s{{\bf 0}}
\def\phi{\varphi}

\newtheorem{theorem}{Theorem}[section]
\newtheorem{lemma1}[theorem]{Lemma}
\newtheorem{remark1}[theorem]{Remark}

\newenvironment{lemma}{\begin{lemma1}}{\end{lemma1}}
\newenvironment{remark}{\begin{remark1}\rm}{\end{remark1}}

\journal{the arXiv}

\begin{document}

\begin{frontmatter}



\title{Fast algorithms for spherical harmonic expansions, III}


\author{Mark Tygert}

\address{Courant Institute of Mathematical Sciences, NYU,
         251 Mercer St., New York, NY 10012}

\begin{abstract}
We accelerate the computation of spherical harmonic transforms,
using what is known as the butterfly scheme.
This provides a convenient alternative to the approach taken
in the second paper from this series
on ``Fast algorithms for spherical harmonic expansions.''
The requisite precomputations become manageable when organized
as a ``depth-first traversal'' of the program's control-flow graph,
rather than as the perhaps more natural ``breadth-first traversal''
that processes one-by-one each level of the multilevel procedure.
We illustrate the results via several numerical examples.
\end{abstract}


\end{frontmatter}



\section{Introduction}

The butterfly algorithm, introduced in~\cite{michielssen-boag}
and~\cite{oneil-woolfe-rokhlin}, is a procedure for rapidly applying
certain matrices to arbitrary vectors.
(Section~\ref{prelims} below provides a brief introduction to the butterfly.)
The present paper uses the butterfly method
in order to accelerate spherical harmonic transforms.
The butterfly procedure does not require the use of extended-precision
arithmetic in order to attain accuracy very close to the machine precision,
not even in its precomputations --- unlike the alternative approach taken
in the predecessor~\cite{tygert_sph} of the present paper.

Unlike some previous works on the butterfly,
the present article does not use on-the-fly evaluation
of individual entries of the matrices
whose applications to vectors are being accelerated.
Instead, we require only efficient evaluation of full columns of the matrices,
in order to make the precomputations affordable.
Furthermore, efficient evaluation of full columns enables
the acceleration of the application to vectors
of both the matrices and their transposes.
On-the-fly evaluation of columns of the matrices
associated with spherical harmonic transforms is available
via the three-term recurrence relations
satisfied by associated Legendre functions
(see, for example, Section~\ref{spharmonics} below).

The precomputations for the butterfly become affordable when organized
as a ``depth-first traversal'' of the program's control-flow graph,
rather than as the perhaps more natural ``breadth-first traversal''
that processes one-by-one each level of the multilevel butterfly procedure
(see Section~\ref{precomps} below).

The present article is supposed to complement~\cite{oneil-woolfe-rokhlin}
and~\cite{tygert_sph}, combining ideas from both.
Although the present paper is self-contained in principle,
we strongly encourage the reader to begin with~\cite{oneil-woolfe-rokhlin}
and~\cite{tygert_sph}.
The original is~\cite{michielssen-boag}.
Major recent developments are in~\cite{candes-demanet-ying} and~\cite{ying}.
The introduction in~\cite{tygert_sph} summarizes
most prior work on computing fast spherical harmonic transforms;
a new application appears in~\cite{reuter-ratner-seideman}.
These articles and their references highlight the computational use
of spherical harmonic transforms in meteorology and quantum chemistry.
The structure of the remainder of the present article is as follows:
Section~\ref{overview} reviews elementary facts
about spherical harmonic transforms.
Section~\ref{prelims} describes basic tools from previous works.
Section~\ref{precomps} organizes the preprocessing for the butterfly
to make memory requirements affordable.
Section~\ref{spharmonics} outlines the application of the butterfly scheme
to the computation of spherical harmonic transforms.
Section~\ref{numerical} describes the results of several numerical tests.
Section~\ref{conclusions} draws some conclusions.

Throughout, we abbreviate ``interpolative decomposition'' to ``ID''
(see Subsection~\ref{ID_subsection} for a description of the ID).
The butterfly procedures formulated in~\cite{michielssen-boag},
\cite{oneil-woolfe-rokhlin}, and the present paper
all use the ID for efficiency.

\section{An overview of spherical harmonic transforms}
\label{overview}

The spherical harmonic expansion of a bandlimited function $f$
on the surface of the sphere has the form
\begin{equation}
\label{expansion}
f(\theta,\phi) = \sum_{k=0}^{2l-1} \sum_{m=-k}^k \beta_k^m
\, \overline{P}_k^{|m|}(\cos(\theta)) \, e^{i m \phi},
\end{equation}
where $(\theta,\phi)$ are the standard spherical coordinates
on the two-dimensional surface of the unit sphere in $\mathbb{R}^3$,
$\theta \in (0,\pi)$ and $\phi \in (0,2\pi)$,
and $\overline{P}_k^{|m|}$ is the normalized associated Legendre function
of degree $k$ and order $|m|$ (see, for example, Subsection~\ref{basic}
for the definition of normalized associated Legendre functions).
Please note that the superscript~$m$ in $\beta^m_k$ denotes an index,
rather than a power. ``Normalized'' refers to the fact that
the normalized associated Legendre functions
of a fixed order are orthonormal on $(-1,1)$
with respect to the standard inner product.
Obviously, the expansion~(\ref{expansion}) contains $4 l^2$ terms.
The complexity of the function $f$ determines $l$.

In many areas of scientific computing, particularly those
using spectral methods for the numerical solution
of partial differential equations,
we need to evaluate the coefficients $\beta^m_k$
in an expansion of the form~(\ref{expansion})
for a function $f$ given by a table of its values
at a collection of appropriately chosen nodes
on the two-dimensional surface of the unit sphere.
Conversely, given the coefficients $\beta^m_k$ in~(\ref{expansion}),
we often need to evaluate $f$ at a collection of points
on the surface of the sphere.
The former is known as the forward spherical harmonic transform,
and the latter is known as the inverse spherical harmonic transform.
A standard discretization of the surface of the sphere
is the ``tensor product,''
consisting of all pairs of the form $(\theta_k,\phi_j)$,
with $\cos(\theta_0)$,~$\cos(\theta_1)$, \dots,
$\cos(\theta_{2l-2})$,~$\cos(\theta_{2l-1})$
being the Gauss-Legendre quadrature nodes of degree $2l$, that is,
\begin{equation}
-1 < \cos(\theta_0) < \cos(\theta_1) < \ldots
   < \cos(\theta_{2l-2}) < \cos(\theta_{2l-1}) < 1
\end{equation}
and
\begin{equation}
\label{thetas}
\overline{P}^0_{2l}(\cos(\theta_k)) = 0
\end{equation}
for $k = 0$,~$1$, \dots, $2l-2$,~$2l-1$,
and with $\phi_0$,~$\phi_1$, \dots, $\phi_{4l-3}$,~$\phi_{4l-2}$
being equispaced on the interval $(0,2\pi)$, that is,
\begin{equation}
\label{phis}
\phi_j = \frac{2\pi \left( j+\frac{1}{2} \right)}{4l-1}
\end{equation}
for $j = 0$,~$1$, \dots, $4l-3$,~$4l-2$.
This leads immediately to numerical schemes
for both the forward and inverse spherical harmonic transforms
whose costs are proportional to $l^3$.

Indeed, given a function $f$ defined on the two-dimensional surface
of the unit sphere by~(\ref{expansion}),
we can rewrite~(\ref{expansion}) in the form
\begin{equation}
\label{alternate_expansion}
f(\theta,\phi) = \sum_{m=-2l+1}^{2l-1} e^{i m \phi} \sum_{k=|m|}^{2l-1}
\beta^m_k \, \overline{P}^{|m|}_k(\cos(\theta)).
\end{equation}
For a fixed value of $\theta$, each of the sums over $k$
in~(\ref{alternate_expansion}) contains no more than $2l$ terms,
and there are $4l-1$ such sums (one for each value of $m$);
since the inverse spherical harmonic transform involves $2l$ values
$\theta_0$,~$\theta_1$, \dots, $\theta_{2l-2}$,~$\theta_{2l-1}$,
the cost of evaluating all sums over $k$ in~(\ref{alternate_expansion})
is proportional to $l^3$. Once all sums over $k$ have been evaluated,
each sum over $m$ may be evaluated for a cost proportional to $l$
(since each of them contains $4l-1$ terms),
and there are $(2l)(4l-1)$ such sums to be evaluated
(one for each pair $(\theta_k,\phi_j)$),
leading to costs proportional to $l^3$ for the evaluation of
all sums over $m$ in~(\ref{alternate_expansion}).
The cost of the evaluation of the whole inverse spherical harmonic transform
(in the form~(\ref{alternate_expansion})) is the sum
of the costs for the sums over $k$ and the sums over $m$,
and is also proportional to $l^3$;
a virtually identical calculation shows that the cost of evaluating
of the forward spherical harmonic transform is also proportional to $l^3$.

A trivial modification of the scheme described in the preceding paragraph
uses the fast Fourier transform (FFT) to evaluate the sums over $m$
in~(\ref{alternate_expansion}),
approximately halving the operation count of the entire procedure.
Several other careful considerations
(see, for example,~\cite{adams-swarztrauber} and~\cite{swarztrauber-spotz})
are able to reduce the costs by 50\% or so,
but there is no simple trick for reducing the costs
of the whole spherical harmonic transform (either forward or inverse)
below $l^3$. The present paper presents faster (albeit more complicated)
algorithms for both forward and inverse spherical harmonic transforms.
Specifically, the present article provides a fast algorithm
for evaluating a sum over $k$ in~(\ref{alternate_expansion})
at $\theta = \theta_0$,~$\theta_1$, \dots, $\theta_{2l-2}$,~$\theta_{2l-1}$,
given the coefficients $\beta^m_{|m|}$,~$\beta^m_{|m|+1}$, \dots,
$\beta^m_{2l-2}$,~$\beta^m_{2l-1}$, for a fixed $m$.
Moreover, the present paper provides a fast algorithm
for the inverse procedure of determining the coefficients
$\beta^m_{|m|}$,~$\beta^m_{|m|+1}$, \dots, $\beta^m_{2l-2}$,~$\beta^m_{2l-1}$
from the values of a sum over $k$ in~(\ref{alternate_expansion})
at $\theta = \theta_0$,~$\theta_1$, \dots, $\theta_{2l-2}$,~$\theta_{2l-1}$.
FFTs or fast discrete sine and cosine transforms can be used
to handle the sums over $m$ in~(\ref{alternate_expansion}) efficiently.
See~\cite{reuter-ratner-seideman} for a detailed summary and novel application
of the overall method. The present article modifies portions of the method
of~\cite{reuter-ratner-seideman} and~\cite{tygert_sph},
focusing exclusively on the modifications.

\section{Preliminaries}
\label{prelims}

In this section, we summarize certain facts from mathematical
and numerical analysis, used in Sections~\ref{precomps} and~\ref{spharmonics}.
Subsection~\ref{ID_subsection} describes interpolative decompositions (IDs).
Subsection~\ref{butterfly} outlines the butterfly algorithm.
Subsection~\ref{basic} summarizes basic properties of normalized
associated Legendre functions.

\subsection{Interpolative decompositions}
\label{ID_subsection}

In this subsection, we define interpolative decompositions (IDs)
and summarize their properties.

The following lemma states that, for any $m \times n$ matrix $A$
of rank $k$, there exist an $m \times k$ matrix $A^{(k)}$
whose columns constitute a subset of the columns of $A$,
and a $k \times n$ matrix $\widetilde{A}$, such that
\begin{enumerate}
\item some subset of the columns of $\widetilde{A}$ makes up
the $k \times k$ identity matrix,
\item $\widetilde{A}$ is not too large, and
\item $A^{(k)}_{m \times k} \cdot \widetilde{A}_{k \times n} = A_{m \times n}$.
\end{enumerate}
Moreover, the lemma provides an approximation 
\begin{equation} \label{2.1}
A^{(k)}_{m \times k} \cdot \widetilde{A}_{k \times n} \approx A_{m \times n}
\end{equation}
when the exact rank of $A$ is greater than $k$,
yet the ($k$+$1$)st greatest singular value of $A$ is still small.
The lemma is a reformulation
of Theorem~3.2 in~\cite{martinsson-rokhlin-tygert1}
and Theorem~3 in~\cite{cheng-gimbutas-martinsson-rokhlin};
its proof is based on techniques described
in~\cite{goreinov-tyrtyshnikov}, \cite{gu-eisenstat96}, 
and~\cite{tyrtyshnikov}.
We will refer to the approximation in~(\ref{2.1}) of $A$
as an interpolative decomposition (ID).
We call $\widetilde{A}$ the ``interpolation matrix'' of the ID.

\begin{lemma}
\label{interpolation_lemma}
Suppose that $m$ and $n$ are positive integers,
and $A$ is a real $m \times n$ matrix.

Then, for any positive integer $k$ with $k \le m$ and $k \le n$,
there exist a real $k \times n$ matrix $\widetilde{A}$,
and a real $m \times k$ matrix $A^{(k)}$ whose columns constitute a subset
of the columns of $A$,
such that
\begin{enumerate}
\item some subset of the columns of $\widetilde{A}$ makes up
the $k \times k$ identity matrix,
\item no entry of $\widetilde{A}$ has an absolute value greater than $1$,
\item the spectral norm (that is, the $l^2$-operator norm) of $\widetilde{A}$
satisfies
$\left\| \widetilde{A}_{k \times n} \right\|_2 \le \sqrt{k (n-k) + 1}$,
\item the least (that is, the $k$\/\,th greatest) singular value
of $\widetilde{A}$ is at least $1$,
\item $A^{(k)}_{m \times k} \cdot \widetilde{A}_{k \times n} = A_{m \times n}$
when $k = m$ or $k = n$, and
\item when $k < m$ and $k < n$, the spectral norm
(that is, the $l^2$-operator norm)
of $A^{(k)}_{m \times k} \cdot \widetilde{A}_{k \times n} - A_{m \times n}$
satisfies
\begin{equation}
\left\| A^{(k)}_{m \times k} \cdot \widetilde{A}_{k \times n}
- A_{m \times n} \right\|_2 \le \sqrt{k (n-k) + 1} \; \sigma_{k+1},
\end{equation}
where $\sigma_{k+1}$ is the ($k$+$1$)st greatest singular value of $A$.
\end{enumerate}
\end{lemma}

Properties~1,~2,~3, and~4 in Lemma~\ref{interpolation_lemma}
ensure that the ID $A^{(k)} \cdot \widetilde{A}$
of $A$ is a numerically stable representation.
Also, property~3 follows directly from properties~1 and~2,
and property~4 follows directly from property~1.

\begin{remark}
\label{QR_algorithm}
Existing algorithms for the computation of the matrices $A^{(k)}$
and $\widetilde{A}$ in Lemma~\ref{interpolation_lemma}
are computationally expensive.
We use instead the algorithm of~\cite{cheng-gimbutas-martinsson-rokhlin}
and~\cite{gu-eisenstat96}
to produce matrices $A^{(k)}$ and $\widetilde{A}$
which satisfy slightly weaker conditions than those
in Lemma~\ref{interpolation_lemma}.
We compute $A^{(k)}$ and $\widetilde{A}$ such that
\begin{enumerate}
\item some subset of the columns of $\widetilde{A}$ makes up
the $k \times k$ identity matrix,
\item \label{computable_bound1}
no entry of $\widetilde{A}$ has an absolute value greater than $2$,
\item the spectral norm (that is, the $l^2$-operator norm) of $\widetilde{A}$
satisfies
$\left\| \widetilde{A}_{k \times n} \right\|_2 \le \sqrt{4 k (n-k) + 1}$,
\item the least (that is, the $k$th greatest) singular value
of $\widetilde{A}$ is at least $1$,
\item \label{exact}
$A^{(k)}_{m \times k} \cdot \widetilde{A}_{k \times n} = A_{m \times n}$
when $k = m$ or $k = n$, and
\item when $k < m$ and $k < n$, the spectral norm
(that is, the $l^2$-operator norm)
of $A^{(k)}_{m \times k} \cdot \widetilde{A}_{k \times n} - A_{m \times n}$
satisfies
\begin{equation}
\label{computable_bound2}
\left\| A^{(k)}_{m \times k} \cdot \widetilde{A}_{k \times n}
- A_{m \times n} \right\|_2 \le \sqrt{4 k (n-k) + 1} \; \sigma_{k+1},
\end{equation}
where $\sigma_{k+1}$ is the ($k$+$1$)st greatest singular value of $A$.
\end{enumerate}
For any positive real number $\epsilon$,
the algorithm can identify the least $k$ such that
$\left\| A^{(k)} \cdot \widetilde{A} - A \right\|_2 \approx \epsilon$.
Furthermore, the algorithm computes both $A^{(k)}$ and $\widetilde{A}$
using at most
\begin{equation}
\label{ID}
C_{\rm ID} = \bigoh(k m n \log(n))
\end{equation}
floating-point operations, typically requiring only
\begin{equation}
\label{ID2}
C'_{\rm ID} = \bigoh(k m n).
\end{equation}
\end{remark}

\subsection{The butterfly algorithm}
\label{butterfly}

In this subsection, we outline a simple case of the butterfly algorithm
from~\cite{michielssen-boag} and~\cite{oneil-woolfe-rokhlin};
see~\cite{oneil-woolfe-rokhlin} for a detailed description.

Suppose that $n$ is a positive integer, and $A$ is an $n \times n$ matrix.
Suppose further that $\epsilon$ and $C$ are positive real numbers,
and $k$ is a positive integer, such that
any contiguous rectangular subblock of $A$ containing at most $Cn$ entries
can be approximated to precision $\epsilon$ by a matrix whose rank is $k$
(using the Frobenius/Hilbert-Schmidt norm to measure the accuracy
of the approximation);
we will refer to this hypothesis as ``{\it the rank property}.''
The running-time of the algorithm will be proportional to $k^2/C$;
taking $C$ to be roughly proportional to $k$ suffices
for many matrices of interest
(including nonequispaced and discrete Fourier transforms),
so ideally $k$ should be small.
We will say that two matrices $G$ and $H$ are equal
to precision $\epsilon$,
denoted $G \approx H$, to mean that the spectral norm
(that is, the $l^2$-operator norm) of $G-H$
is $\bigoh(\epsilon)$.

We now explicitly use the rank property for subblocks of multiple heights,
to illustrate the basic structure of the butterfly scheme.

Consider any two adjacent contiguous rectangular subblocks $L$ and $R$ of $A$,
each containing at most $Cn$ entries and having the same numbers of rows,
with $L$ on the left and $R$ on the right.
Due to the rank property, there exist IDs
\begin{equation}
\label{1}
L \approx L^{(k)} \cdot \widetilde{L}
\end{equation}
and
\begin{equation}
\label{2}
R \approx R^{(k)} \cdot \widetilde{R},
\end{equation}
where $L^{(k)}$ is a matrix having $k$ columns,
which constitute a subset of the columns of $L$,
$R^{(k)}$ is a matrix having $k$ columns,
which constitute a subset of the columns of $R$,
$\widetilde{L}$ and $\widetilde{R}$ are matrices each having $k$ rows,
and all entries of $\widetilde{L}$ and $\widetilde{R}$ have absolute values
of at most 2.

To set notation, we concatenate the matrices $L$ and $R$,
and split the columns of the result in half (or approximately in half),
obtaining $T$ on top and $B$ on the bottom:
\begin{equation}
\label{3}
\left( \begin{array}{c|c} L & R \end{array} \right)
= \left( \begin{array}{c} T \\\hline B \end{array} \right).
\end{equation}
Observe that the matrices $T$ and $B$ each have at most $Cn$ entries
(since $L$ and $R$ each have at most $Cn$ entries).
Similarly, we concatenate the matrices $L^{(k)}$ and $R^{(k)}$,
and split the columns of the result in half (or approximately in half),
obtaining $T^{(2k)}$ and $B^{(2k)}$:
\begin{equation}
\label{4}
\left( \begin{array}{c|c} L^{(k)} & R^{(k)} \end{array} \right)
= \left( \begin{array}{c} T^{(2k)} \\\hline B^{(2k)} \end{array} \right).
\end{equation}
Observe that the $2k$ columns of $T^{(2k)}$ are also columns of $T$,
and that the $2k$ columns of $B^{(2k)}$ are also columns of $B$.

Due to the rank property, there exist IDs
\begin{equation}
\label{5}
T^{(2k)} \approx T^{(k)} \cdot \widetilde{T^{(2k)}}
\end{equation}
and
\begin{equation}
\label{6}
B^{(2k)} \approx B^{(k)} \cdot \widetilde{B^{(2k)}},
\end{equation}
where $T^{(k)}$ is a matrix having $k$ columns,
which constitute a subset of the columns of $T^{(2k)}$,
$B^{(k)}$ is a matrix having $k$ columns,
which constitute a subset of the columns of $B^{(2k)}$,
$\widetilde{T^{(2k)}}$ and $\widetilde{B^{(2k)}}$ are matrices
each having $k$ rows, and all entries of $\widetilde{T^{(2k)}}$
and $\widetilde{B^{(2k)}}$ have absolute values of at most 2.

Combining~(\ref{1})--(\ref{6}) yields that
\begin{equation}
\label{top}
T \approx T^{(k)} \cdot \widetilde{T^{(2k)}}
\cdot \left( \begin{array}{c|c} \widetilde{L} & \0s \\\hline
                                \0s & \widetilde{R} \end{array} \right)
\end{equation}
and
\begin{equation}
\label{bottom}
B \approx B^{(k)} \cdot \widetilde{B^{(2k)}}
\cdot \left( \begin{array}{c|c} \widetilde{L} & \0s \\\hline
                                \0s & \widetilde{R} \end{array} \right).
\end{equation}

If we use $m$ to denote the number of rows in $L$ (which is the same
as the number of rows in $R$),
then the number of columns in $L$ (or $R$) is at most $Cn/m$,
and so the total number of entries in the matrices
in the right-hand sides of~(\ref{1}) and~(\ref{2})
can be as large as $2mk + 2k(Cn/m)$,
whereas the total number of nonzero entries in the matrices
in the right-hand sides of~(\ref{top}) and~(\ref{bottom})
is at most $mk + 4k^2 + 2k(Cn/m)$.
If $m$ is nearly as large as possible --- nearly $n$ ---
and $k$ and $C$ are much smaller than $n$,
then $mk + 4k^2 + 2k(Cn/m)$ is about half $2mk + 2k(Cn/m)$.
Thus, the representation provided in~(\ref{top}) and~(\ref{bottom})
of the merged matrix from~(\ref{3})
is more efficient than that provided in~(\ref{1}) and~(\ref{2}),
both in terms of the memory required for storage,
and in terms of the number of operations required for applications to vectors.
Notice the advantage of using the rank property for blocks of multiple heights.

\begin{figure}
\begin{center}
\caption{The partitionings in the multilevel decomposition
         for an $8 \times 8$ matrix with $C = 1$}
\label{multi}
\vspace{1.5em}
\begin{minipage}{1in}
\begin{center}
\scalebox{.25}{\includegraphics{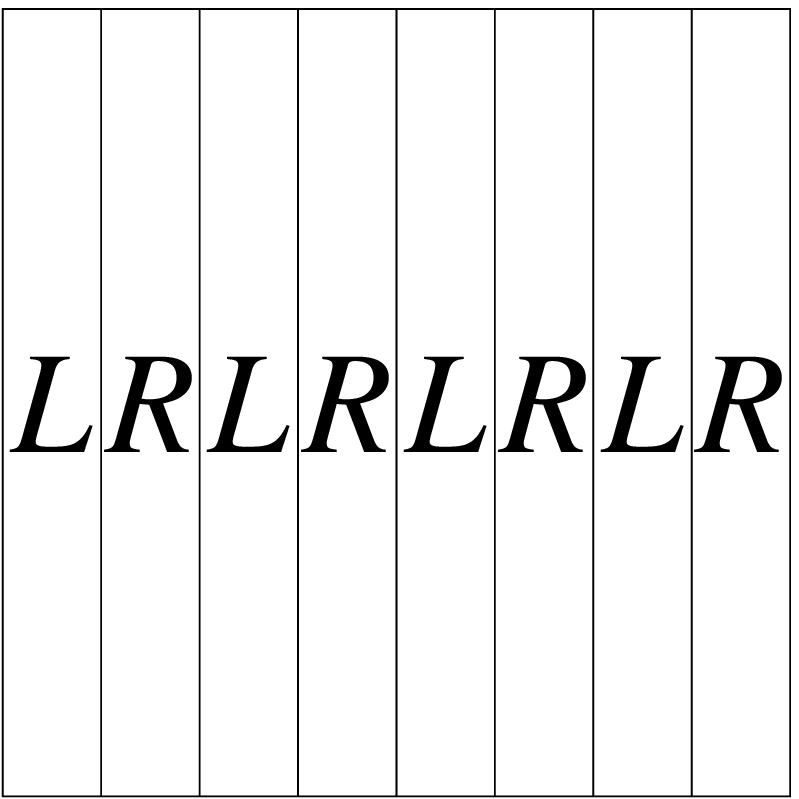}} \\
Level 1
\end{center}
\end{minipage}
\begin{minipage}{1in}
\begin{center}
\scalebox{.25}{\includegraphics{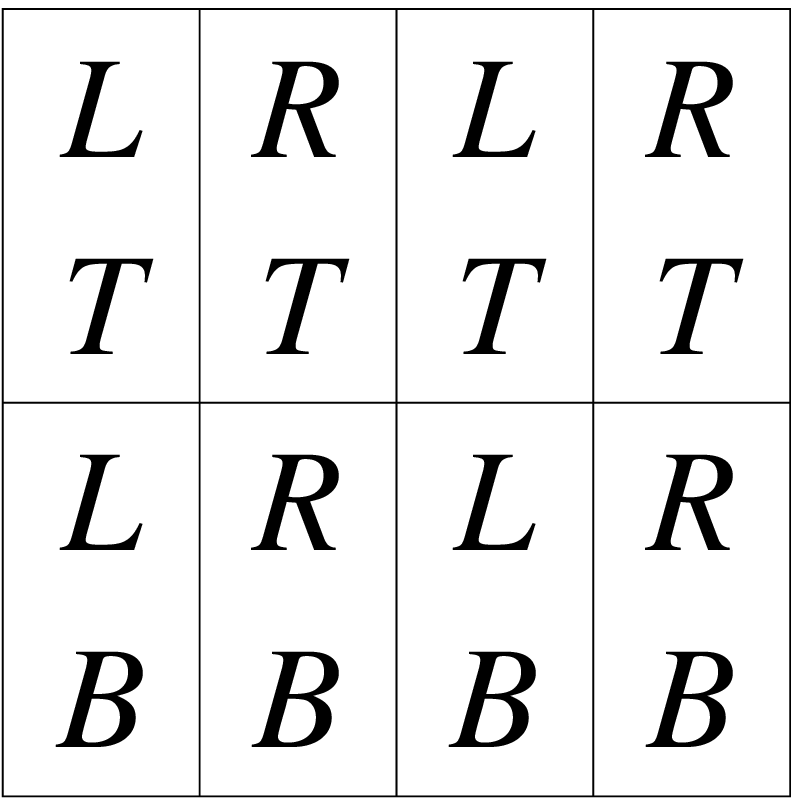}} \\
Level 2
\end{center}
\end{minipage}
\begin{minipage}{1in}
\begin{center}
\scalebox{.25}{\includegraphics{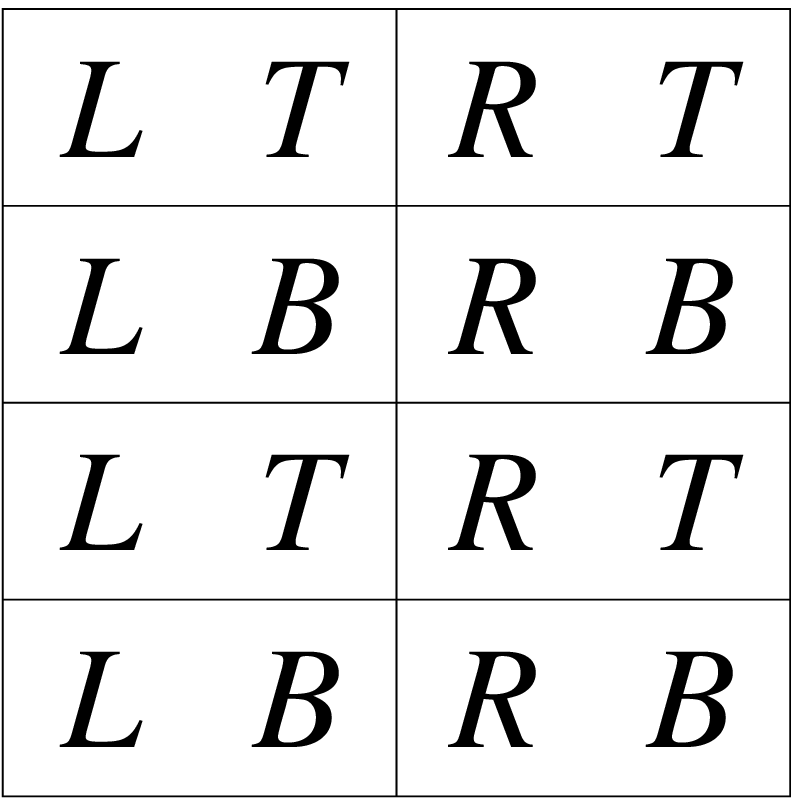}} \\
Level 3
\end{center}
\end{minipage}
\begin{minipage}{1in}
\begin{center}
\scalebox{.25}{\includegraphics{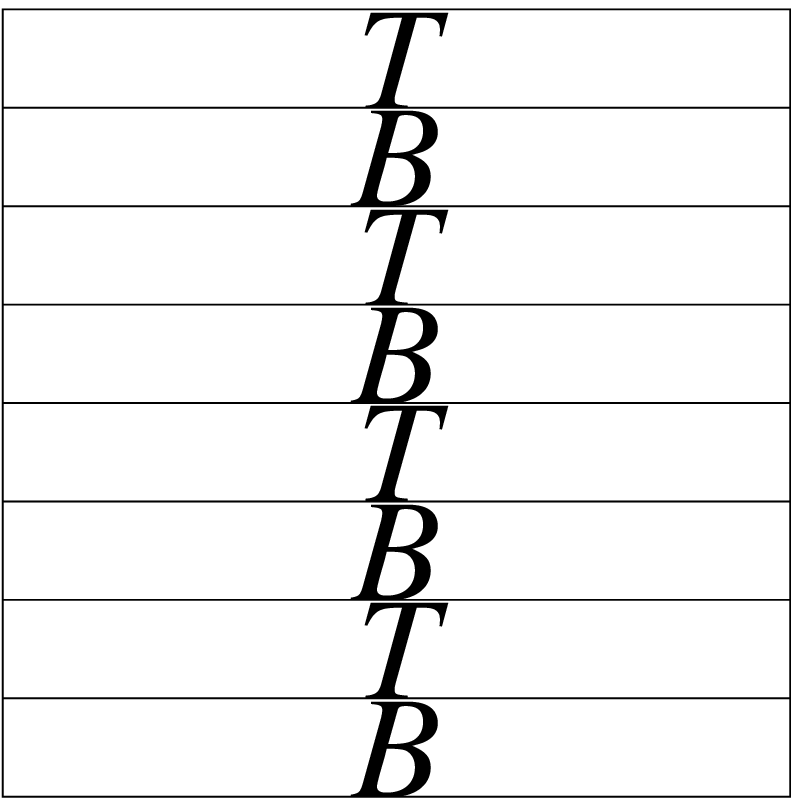}} \\
Level 4
\end{center}
\end{minipage}
\\\vspace{1 em}
$L$ indicates the left member of a pair; $R$ indicates the right member. \\
$T$ indicates the top member of a pair; $B$ indicates the bottom member.
\end{center}
\end{figure}

Naturally, we may repeat this process of merging adjacent blocks
and splitting in half the columns of the result,
updating the compressed representations after every split.
We start by partitioning $A$ into blocks
each dimensioned $n \times \lfloor C \rfloor$
(except possibly for the rightmost block, which may have fewer 
than $\lfloor C \rfloor$ columns),
and then repeatedly group unprocessed blocks (of whatever dimensions)
into disjoint pairs, processing these pairs by merging and splitting them
into new, unprocessed blocks having fewer rows.
The resulting multilevel representation of $A$ allows us to apply $A$
with precision $\epsilon$
from the left to any column vector, or from the right to any row vector,
using just $\bigoh((k^2/C) \, n \log(n))$ floating-point operations
(there will be $\bigoh(\log(n))$ levels in the scheme,
and each level except for the last
will only involve $\bigoh(n/C)$ interpolation matrices
of dimensions $k \times (2k)$,
such as $\widetilde{T^{(2k)}}$ and $\widetilde{B^{(2k)}}$).
Figure~\ref{multi} illustrates the resulting partitionings of $A$ into blocks
of various dimensions (but with every block having the same number of entries),
when $n = 8$ and $C = 1$.
For further details, see~\cite{oneil-woolfe-rokhlin}.

\begin{remark}
Needless to say, the same multilevel representation of $A$
permits the rapid application of $A$ both from the left to column vectors
and from the right to row vectors.
There is no need for constructing multilevel representations
of both $A$ and the transpose of $A$.
\end{remark}

\begin{remark}
\label{adaptivity}
In practice, the IDs used for accurately approximating subblocks of $A$
do not all have the same fixed rank~$k$. Instead, for each subblock,
we determine the minimal possible rank such that the associated ID still
approximates the subblock to precision $\epsilon$, and we use this ID
in place of one whose rank is $k$.
Determining ranks adaptively in this manner accelerates the algorithm
substantially. For further details, see~\cite{oneil-woolfe-rokhlin}.
All our implementations use this adaptation.
\end{remark}

\subsection{Normalized associated Legendre functions}
\label{basic}

In this subsection, we discuss several classical facts concerning
normalized associated Legendre functions.
All of these facts follow trivially from results contained,
for example, in~\cite{abramowitz-stegun} or~\cite{szego}.

For any nonnegative integers $l$ and $m$ such that $l \ge m$,
we use $\overline{P}_l^m$ to denote the normalized associated Legendre function
of degree $l$ and order $m$, defined on $(-1,1)$ via the formula
\begin{equation}
\label{association}
\overline{P}^m_l(x) = \sqrt{\frac{2l+1}{2} \frac{(l-m)!}{(l+m)!}}
                   \; \left(1-x^2\right)^{m/2} \; \frac{d^m}{dx^m}P_l(x),
\end{equation}
where $P_l$ is the Legendre polynomial of degree $l$,
\begin{equation}
P_l(x) = \frac{1}{2^l \, l!} \; \frac{d^l}{dx^l}\left((x^2-1)^l\right)
\end{equation}
(see, for example, Chapter~8 of~\cite{abramowitz-stegun}).
``Normalized'' refers to the fact that
the normalized associated Legendre functions
of a fixed order $m$ are orthonormal on $(-1,1)$
with respect to the standard inner product.
If $l-m$ is even, then $\overline{P}^m_l(-x) = \overline{P}^m_l(x)$
for any $x \in (-1,1)$.
If $l-m$ is odd, then $\overline{P}^m_l(-x) = -\overline{P}^m_l(x)$
for any $x \in (-1,1)$.

The following lemma states that the normalized associated Legendre functions
satisfy a certain self-adjoint second-order linear (Sturm-Liouville)
differential equation.

\begin{lemma}
Suppose that $m$ is a nonnegative integer.

Then,
\begin{equation}
\label{sturm-liouville_Leg}
-\frac{d}{dx}\left( (1-x^2) \, \frac{d}{dx}\overline{P}^m_l(x) \right)
+\left( \frac{m^2}{1-x^2} - l(l+1) \right) \, \overline{P}^m_l(x) = 0
\end{equation}
for any $x \in (-1,1)$, and $l = m$,~$m+1$, $m+2$, \dots.
\end{lemma}

The following lemma states that the normalized associated Legendre function
of order $m$ and degree $m+2n$ has exactly $n$ zeros inside $(0,1)$,
and, moreover, that the normalized associated Legendre function
of order $m$ and degree $m+2n+1$ also has exactly $n$ zeros inside $(0,1)$.

\begin{lemma}
\label{distinct_theorem}
Suppose that $m$ and $n$ are nonnegative integers with $n > 0$.

Then, there exist precisely $n$ real numbers
$x_0$,~$x_1$, \dots, $x_{n-2}$,~$x_{n-1}$ such that
\begin{equation}
\label{distinctness}
0 < x_0 < x_1 < \ldots < x_{n-2} < x_{n-1} < 1
\end{equation}
and
\begin{equation}
\label{zeros}
\overline{P}^m_{m+2n}(x_j) = 0
\end{equation}
for $j = 0$,~$1$, \dots, $n-2$,~$n-1$.

Moreover, there exist precisely $n$ real numbers
$y_0$,~$y_1$, \dots, $y_{n-2}$,~$y_{n-1}$ such that
\begin{equation}
\label{distinctness2}
0 < y_0 < y_1 < \ldots < y_{n-2} < y_{n-1} < 1
\end{equation}
and
\begin{equation}
\label{zeros2}
\overline{P}^m_{m+2n+1}(y_j) = 0
\end{equation}
for $j = 0$,~$1$, \dots, $n-2$,~$n-1$.
\end{lemma}

Suppose that $m$ and $n$ are nonnegative integers with $n > 0$.
Then, we define real numbers
$\rho_0$,~$\rho_1$, \dots, $\rho_{n-2}$,~$\rho_{n-1}$,
$\sigma_0$,~$\sigma_1$, \dots, $\sigma_{n-2}$,~$\sigma_{n-1}$,
and $\sigma_n$ via the formulae
\begin{equation}
\label{weights1}
\rho_j
= \frac{2 \, (2m+4n+1)}
       {\left(1-(x_j)^2\right)
        \, \left( \frac{d}{dx} \overline{P}^m_{m+2n}(x_j) \right)^2}
\end{equation}
for $j = 0$,~$1$, \dots, $n-2$,~$n-1$,
where $x_0$,~$x_1$, \dots, $x_{n-2}$,~$x_{n-1}$ are from~(\ref{zeros}),
\begin{equation}
\label{weights2}
\sigma_j
= \frac{2 \, (2m+4n+3)}
  {\left(1-(y_j)^2\right)
   \, \left( \frac{d}{dx} \overline{P}^m_{m+2n+1}(y_j) \right)^2}
\end{equation}
for $j = 0$,~$1$, \dots, $n-2$,~$n-1$,
where $y_0$,~$y_1$, \dots, $y_{n-2}$,~$y_{n-1}$ are from~(\ref{zeros2}), and
\begin{equation}
\label{weights2_0}
\sigma_n
= \frac{2m+4n+3}{\left( \frac{d}{dx} \overline{P}^m_{m+2n+1}(0) \right)^2}.
\end{equation}

The following lemma describes what are known as Gauss-Jacobi quad\-rature
formulae corresponding to associated Legendre functions.

\begin{lemma}
\label{quadratures}
Suppose that $m$ and $n$ are nonnegative integers with $n > 0$.

Then,
\begin{equation}
\label{Christoffel_definition}
\int_{-1}^1 dx \; \left(1-x^2\right)^m \; p(x)
= \sum_{j=0}^{n-1} \rho_j \; \left(1-(x_j)^2\right)^m \; p(x_j)
\end{equation}
for any even polynomial $p$ of degree at most $4n-2$,
where $x_0$,~$x_1$, \dots, $x_{n-2}$,~$x_{n-1}$ are from~(\ref{zeros}),
and $\rho_0$,~$\rho_1$, \dots, $\rho_{n-2}$,~$\rho_{n-1}$
are defined in~(\ref{weights1}).

Furthermore,
\begin{equation}
\label{Christoffel_definition2}
\int_{-1}^1 dx \; \left(1-x^2\right)^m \; p(x)
= \sigma_n \; p(0)
+ \sum_{j=0}^{n-1} \sigma_j \; \left(1-(y_j)^2\right)^m \; p(y_j)
\end{equation}
for any even polynomial $p$ of degree at most $4n$,
where $y_0$,~$y_1$, \dots, $y_{n-2}$,~$y_{n-1}$ are from~(\ref{zeros2}),
and $\sigma_0$,~$\sigma_1$, \dots, $\sigma_{n-1}$,~$\sigma_n$
are defined in~(\ref{weights2}) and~(\ref{weights2_0}).
\end{lemma}

\begin{remark}
Formulae~(35) and~(36) of~\cite{tygert_sph} incorrectly
omitted the factors $(1-(x_j)^2)^m$ and $(1-(y_j)^2)^m$
appearing in the analogous~(\ref{Christoffel_definition})
and~(\ref{Christoffel_definition2}) above.
\end{remark}

Suppose that $m$ is a nonnegative integer.
Then, we define real numbers $c_m$,~$c_{m+1}$, $c_{m+2}$, \dots\
and $d_m$,~$d_{m+1}$, $d_{m+2}$, \dots\ via the formulae
\begin{equation}
\label{c_rec}
c_l = \sqrt{\frac{(l-m+1)(l-m+2)(l+m+1)(l+m+2)}
           {(2l+1) \, (2l+3)^2 \, (2l+5)}}
\end{equation}
for $l = m$,~$m+1$,~$m+2$, \dots, and
\begin{equation}
\label{d_rec}
d_l = \frac{2 l (l+1) - 2 m^2 - 1}{(2l-1)(2l+3)}
\end{equation}
for $l = m$,~$m+1$,~$m+2$, \dots.

The following lemma states that the normalized associated Legendre functions
of a fixed order $m$ satisfy a certain three-term recurrence relation.

\begin{lemma}
\label{rec_relation}
Suppose that $m$ is a nonnegative integer.

Then,
\begin{equation}
\label{recurrence0}
x^2 \; \overline{P}^m_l(x)
= d_l \; \overline{P}^m_l(x)
+ c_l \; \overline{P}^m_{l+2}(x)
\end{equation}
for any $x \in (-1,1)$, and $l = m$ or $l = m+1$, and
\begin{equation}
\label{recurrence}
x^2 \; \overline{P}^m_l(x)
= c_{l-2} \; \overline{P}^m_{l-2}(x)
+ d_l \; \overline{P}^m_l(x)
+ c_l \; \overline{P}^m_{l+2}(x)
\end{equation}
for any $x \in (-1,1)$, and $l = m+2$,~$m+3$,~$m+4$, \dots,
where $c_m$,~$c_{m+1}$,~$c_{m+2}$, \dots\ are defined in~(\ref{c_rec}),
and $d_m$,~$d_{m+1}$,~$d_{m+2}$, \dots\ are defined in~(\ref{d_rec}).
\end{lemma}

\section{Precomputations for the butterfly scheme}
\label{precomps}

In this section, we discuss the preprocessing required
for the butterfly algorithm summarized in Subsection~\ref{butterfly}.
We will be using the notation detailed in Subsection~\ref{butterfly}.

Perhaps the most natural organization of the computations required
to construct the multilevel representation of an $n \times n$ matrix $A$
is first to process all blocks having $n$ rows
(Level~1 in Figure~\ref{multi} above),
then to process all blocks having about $n/2$ rows
(Level~2 in Figure~\ref{multi}),
then to process all blocks having about $n/4$ rows
(Level~3 in Figure~\ref{multi}), and so on.
Indeed, \cite{oneil-woolfe-rokhlin} uses this organization,
which amounts to a ``breadth-first traversal''
of the control-flow graph for the program applying $A$ to a vector
(see, for example, \cite{aho-hopcroft-ullman}
for an introduction to ``breadth-first'' and ``depth-first'' orderings).
This scheme for preprocessing is efficient when the entries of $A$
can be efficiently computed on-the-fly, individually.
(Of course, we are assuming that $A$ has a suitable rank property, that is,
that there are positive real numbers $\epsilon$ and $C$,
and a positive integer $k$, such that
any contiguous rectangular subblock of $A$ containing at most $Cn$ entries
can be approximated to precision $\epsilon$ by a matrix whose rank is $k$,
using the Frobenius/Hilbert-Schmidt norm to measure the accuracy
of the approximation. Often, taking $C$ to be roughly proportional
to $k$ suffices, and ideally $k$ and $\epsilon$ are small.)
If the entries of $A$ cannot be efficiently computed individually,
however, then the ``breadth-first traversal'' may need
to store $\bigoh(n^2)$ entries at some point during the precomputations,
in order to avoid recomputing entries of the matrix.

If individual columns of $A$
(but not necessarily arbitrary individual entries) can be computed efficiently,
then ``depth-first traversal'' of the control-flow graph
requires only $\bigoh((k^2/C) \, n \log(n))$ floating-point words of memory
at any point during the precomputations, for the following reason.
We will say that we ``process'' a block of $A$ to mean that we merge it
with another, and split and recompress the result,
producing a pair of new, unprocessed blocks.
Rather than starting the preprocessing by constructing all blocks
having $n$ rows, we construct each such block only after processing
as many blocks as possible which previous processing creates,
but which have not yet been processed.
Furthermore, we construct each block having $n$ rows only after having already
constructed (and possibly processed) all blocks to its left.
To reiterate, we construct a block having $n$ rows only after having exhausted
all possibilities for both creating and processing blocks to its left.

For each {\it processed} block $B$, we need only store the interpolation matrix
$\widetilde{B^{(2k)}}$ and the indices of the columns chosen for the ID;
we need not store the $k$ columns of $B^{(k)}$ selected for the ID,
since the algorithm for applying $A$ (or its transpose) to a vector
never explicitly uses any columns of a block that has been merged
with another and split, but instead interpolates from (or anterpolates to)
the shorter blocks arising from the processing.
Conveniently, the matrix $\widetilde{B^{(2k)}}$ that we must store is small
-- no larger than $k \times (2k)$.
For each {\it unprocessed} block $B$, we do need to store the $k$ columns
in $B^{(k)}$ selected for the ID, in addition to storing $\widetilde{B^{(2k)}}$
and the indices of the columns chosen for the ID,
facilitating any subsequent processing.
Although $B^{(k)}$ may have many rows, it has only as many rows as $B$
and hence is smaller when $B$ has fewer rows.
Thus, every time we process a pair of tall blocks, producing a new pair
of blocks having half as many rows,
the storage requirements for all these blocks together nearly halve.
By always processing as many already constructed blocks as possible,
we minimize the amount of memory required.

\section{Spherical harmonic transforms via the butterfly scheme}
\label{spharmonics}

In this section, we describe how to use the butterfly algorithm
to compute fast spherical harmonic transforms,
via appropriate modifications of the algorithm of~\cite{tygert_sph}.

We substitute the butterfly algorithm for the divide-and-conquer algorithm
of~\cite{gu-eisenstat95} used in Section~3.1 of~\cite{tygert_sph},
otherwise leaving the approach of~\cite{tygert_sph} unchanged. Specifically,
given numbers $\beta_0$,~$\beta_1$, \dots, $\beta_{n-2}$,~$\beta_{n-1}$,
we use the butterfly scheme to compute the numbers
$\alpha_0$,~$\alpha_1$, \dots, $\alpha_{n-2}$,~$\alpha_{n-1}$ defined
via the formula
\begin{equation}
\label{transform}
\alpha_i
= \sum_{j=0}^{n-1} \beta_j \, \sqrt{\rho_i} \;\; \overline{P}^m_{m+2j}(x_i),
\end{equation}
for $i = 0$,~$1$, \dots, $n-2$,~$n-1$,
where $m$ is a nonnegative integer,
$\overline{P}^m_m$,~$\overline{P}^m_{m+2}$, \dots,
$\overline{P}^m_{m+2n-2}$,~$\overline{P}^m_{m+2n}$
are the normalized associated Legendre functions of order $m$ defined
in~(\ref{association}),
$x_0$,~$x_1$, \dots, $x_{n-2}$,~$x_{n-1}$ are the positive zeros
of $\overline{P}^m_{m+2n}$ from~(\ref{zeros}),
and $\rho_0$,~$\rho_1$, \dots, $\rho_{n-2}$,~$\rho_{n-1}$
are the corresponding quadrature weights from~(\ref{Christoffel_definition}).
Similarly,
given numbers~$\alpha_0$,~$\alpha_1$, \dots, $\alpha_{n-2}$,~$\alpha_{n-1}$,
we use the butterfly scheme to compute the numbers
$\beta_0$,~$\beta_1$, \dots, $\beta_{n-2}$,~$\beta_{n-1}$
satisfying~(\ref{transform}).
The factors $\sqrt{\rho_0}$,~$\sqrt{\rho_1}$, \dots,
$\sqrt{\rho_{n-2}}$,~$\sqrt{\rho_{n-1}}$ ensure that the linear transformation
mapping $\beta_0$,~$\beta_1$, \dots, $\beta_{n-2}$,~$\beta_{n-1}$
to $\alpha_0$,~$\alpha_1$, \dots, $\alpha_{n-2}$,~$\alpha_{n-1}$
via~(\ref{transform}) is unitary
(due to~(\ref{association}),~(\ref{Christoffel_definition}),
and the orthonormality of the normalized associated Legendre functions
on $(-1,1)$), so that the inverse of the linear transformation
is its transpose.

Moreover, given numbers
$\nu_0$,~$\nu_1$, \dots, $\nu_{n-2}$,~$\nu_{n-1}$,
we use the butterfly scheme to compute the numbers
$\mu_0$,~$\mu_1$, \dots, $\mu_{n-2}$,~$\mu_{n-1}$ defined via the formula
\begin{equation}
\label{transform2}
\mu_i
= \sum_{j=0}^{n-1} \nu_j \, \sqrt{\sigma_i} \;\; \overline{P}^m_{m+2j+1}(y_i),
\end{equation}
for $i = 0$,~$1$, \dots, $n-2$,~$n-1$,
where $m$ is a nonnegative integer,
$\overline{P}^m_{m+1}$,~$\overline{P}^m_{m+3}$, \dots,
$\overline{P}^m_{m+2n-1}$,~$\overline{P}^m_{m+2n+1}$
are the normalized associated Legendre functions of order $m$ defined
in~(\ref{association}),
$y_0$,~$y_1$, \dots, $y_{n-2}$,~$y_{n-1}$ are the positive zeros
of $\overline{P}^m_{m+2n+1}$ from~(\ref{zeros2}),
and $\sigma_0$,~$\sigma_1$, \dots, $\sigma_{n-2}$,~$\sigma_{n-1}$
are the corresponding quadrature weights from~(\ref{Christoffel_definition2}).
Similarly,
given numbers~$\mu_0$,~$\mu_1$, \dots, $\mu_{n-2}$,~$\mu_{n-1}$,
we use the butterfly scheme to compute the numbers
$\nu_0$,~$\nu_1$, \dots, $\nu_{n-2}$,~$\nu_{n-1}$
satisfying~(\ref{transform2}).
As above, the factors $\sqrt{\sigma_0}$,~$\sqrt{\sigma_1}$, \dots,
$\sqrt{\sigma_{n-2}}$,~$\sqrt{\sigma_{n-1}}$ ensure that
the linear transformation mapping
$\nu_0$,~$\nu_1$, \dots, $\nu_{n-2}$,~$\nu_{n-1}$
to $\mu_0$,~$\mu_1$, \dots, $\mu_{n-2}$,~$\mu_{n-1}$
via~(\ref{transform2}) is unitary, so that its inverse is its transpose.

Computing spherical harmonic transforms requires
several additional computations, detailed in~\cite{tygert_sph}.
(See also Remark~\ref{omission} below.)
The butterfly algorithm replaces only the procedure described
in Section~3.1 of~\cite{tygert_sph}.

In order to use~(\ref{transform}) and~(\ref{transform2}) numerically,
we need to precompute the positive zeros
$x_0$,~$x_1$, \dots, $x_{n-2}$,~$x_{n-1}$ 
of $\overline{P}^m_{m+2n}$ from~(\ref{zeros}),
the corresponding quadrature weights
$\rho_0$,~$\rho_1$, \dots, $\rho_{n-2}$,~$\rho_{n-1}$
from~(\ref{Christoffel_definition}),
the positive zeros $y_0$,~$y_1$, \dots, $y_{n-2}$,~$y_{n-1}$
of $\overline{P}^m_{m+2n+1}$ from~(\ref{zeros2}),
and the corresponding quadrature weights
$\sigma_0$,~$\sigma_1$, \dots, $\sigma_{n-2}$,~$\sigma_{n-1}$
from~(\ref{Christoffel_definition2}).
Section~3.3 of~\cite{tygert_sph} describes suitable procedures
(based on integrating the ordinary differential equation
in~(\ref{sturm-liouville_Leg}) in ``Pr\"ufer coordinates'').
We found it expedient to perform this preprocessing in extended-precision
arithmetic, in order to compensate for the loss of a couple of digits
of accuracy relative to the machine precision.

To perform the precomputations described in Section~\ref{precomps} above
associated with~(\ref{transform}) and~(\ref{transform2}),
we need to be able to evaluate efficiently all $n$ functions
$\overline{P}^m_m$,~$\overline{P}^m_{m+2}$, \dots,
$\overline{P}^m_{m+2n-4}$,~$\overline{P}^m_{m+2n-2}$
at any of the precomputed positive zeros
$x_0$,~$x_1$, \dots, $x_{n-2}$,~$x_{n-1}$ 
of $\overline{P}^m_{m+2n}$ from~(\ref{zeros}),
and, similarly, we need to be able to evaluate efficiently all $n$ functions
$\overline{P}^m_{m+1}$,~$\overline{P}^m_{m+3}$, \dots,
$\overline{P}^m_{m+2n-3}$,~$\overline{P}^m_{m+2n-1}$
at any of the precomputed positive zeros
$y_0$,~$y_1$, \dots, $y_{n-2}$,~$y_{n-1}$ 
of $\overline{P}^m_{m+2n+1}$ from~(\ref{zeros2}).
For this, we may use the recurrence relations~(\ref{recurrence0})
and~(\ref{recurrence}), starting with the values of $\overline{P}^m_m(x_i)$,
$\overline{P}^m_{m+1}(y_i)$, $\overline{P}^m_{m+2}(x_i)$,
and~$\overline{P}^m_{m+3}(y_i)$ obtained via~(\ref{association}).
(We can counter underflow by tracking exponents explicitly,
in the standard fashion.)
Such use of the recurrence is a classic procedure;
see, for example, Chapter~8 of~\cite{abramowitz-stegun}.
The recurrence appears to be numerically stable when used
for evaluating normalized associated Legendre functions of order $m$
and of degrees at most $m+2n-1$,
at these special points $x_0$,~$x_1$, \dots, $x_{n-2}$,~$x_{n-1}$ 
and $y_0$,~$y_1$, \dots, $y_{n-2}$,~$y_{n-1}$, even when $n$ is very large.
We did not need to use extended-precision arithmetic for this preprocessing.

\begin{remark}
\label{omission}
The formula~(88) in~\cite{tygert_sph} that is analogous to~(\ref{transform})
of the present paper omits the factors $\sqrt{\rho_0}$,~$\sqrt{\rho_1}$, \dots,
$\sqrt{\rho_{n-2}}$,~$\sqrt{\rho_{n-1}}$ included in~(\ref{transform}).
Obviously, the vectors
\begin{equation}
\label{no_diag}
(\alpha_0, \alpha_1, \dots, \alpha_{n-2}, \alpha_{n-1})^\top
\end{equation}
and
\begin{equation}
\label{with_diag}
\left(\frac{\alpha_0}{\sqrt{\rho_0}}, \frac{\alpha_1}{\sqrt{\rho_1}}, \dots,
\frac{\alpha_{n-2}}{\sqrt{\rho_{n-2}}}, \frac{\alpha_{n-1}}{\sqrt{\rho_{n-1}}}
\right)^\top
\end{equation}
differ by a diagonal transformation, and so we can obtain either one
from the other efficiently.
In fact, the well-conditioned matrix $A$ from Section~3.1 of~\cite{tygert_sph}
represents the same diagonal transformation,
mapping~(\ref{no_diag}) to~(\ref{with_diag}).
Similar remarks apply to~(\ref{transform2}), of course.
\end{remark}

\section{Numerical results}
\label{numerical}

In this section, we describe the results of several numerical tests
of the algorithm of the present paper.
(Computing spherical harmonic transforms requires
several additional computations, detailed in~\cite{tygert_sph}
--- see Section~\ref{spharmonics} for further information.
The butterfly algorithm replaces only the procedure described
in Section~3.1 of~\cite{tygert_sph}.)

Tables~\ref{evens_time}--\ref{small_pre}
report the results of computing from real numbers
$\beta_0$,~$\beta_1$, \dots, $\beta_{n-2}$,~$\beta_{n-1}$
the real numbers $\alpha_0$,~$\alpha_1$, \dots, $\alpha_{n-2}$,~$\alpha_{n-1}$
defined by the formula
\begin{equation}
\label{evens}
\alpha_i = \sum_{j=0}^{n-1} \beta_j \, \sqrt{\rho_i} \;\;
           \overline{P}^m_{m+2j}(x_i)
\end{equation}
for $i = 0$,~$1$, \dots, $n-2$,~$n-1$,
where $\overline{P}^m_m$,~$\overline{P}^m_{m+2}$, \dots,
$\overline{P}^m_{m+2n-2}$,~$\overline{P}^m_{m+2n}$
are the normalized associated Legendre functions defined
in~(\ref{association}),
$x_0$,~$x_1$, \dots, $x_{n-2}$,~$x_{n-1}$ are the positive zeros
of $\overline{P}^m_{m+2n}$ from~(\ref{zeros}),
and $\rho_0$,~$\rho_1$, \dots, $\rho_{n-2}$,~$\rho_{n-1}$
are the corresponding quadrature weights from~(\ref{Christoffel_definition}).
We will refer to the map via~(\ref{evens})
from $\beta_0$,~$\beta_1$, \dots, $\beta_{n-2}$,~$\beta_{n-1}$
to $\alpha_0$,~$\alpha_1$, \dots, $\alpha_{n-2}$,~$\alpha_{n-1}$
as the forward transform, and the map via~(\ref{evens})
from $\alpha_0$,~$\alpha_1$, \dots, $\alpha_{n-2}$,~$\alpha_{n-1}$
to $\beta_0$,~$\beta_1$, \dots, $\beta_{n-2}$,~$\beta_{n-1}$
as the inverse transform (the inverse is also the transpose,
due to~(\ref{association}),~(\ref{Christoffel_definition}),
and the orthonormality of the normalized associated Legendre functions
on $(-1,1)$).
The values of $m$ differ in Tables~1--2 and 3--4.

Tables~\ref{odds_time} and~\ref{odds_pre} report the results
of computing from real numbers
$\nu_0$,~$\nu_1$, \dots, $\nu_{n-2}$,~$\nu_{n-1}$
the real numbers $\mu_0$,~$\mu_1$, \dots, $\mu_{n-2}$,~$\mu_{n-1}$
defined by the formula
\begin{equation}
\label{odds}
\mu_i = \sum_{j=0}^{n-1} \nu_j \, \sqrt{\sigma_i} \;\;
        \overline{P}^m_{m+2j+1}(y_i)
\end{equation}
for $i = 0$,~$1$, \dots, $n-2$,~$n-1$,
where $\overline{P}^m_{m+1}$,~$\overline{P}^m_{m+3}$, \dots,
$\overline{P}^m_{m+2n-1}$,~$\overline{P}^m_{m+2n+1}$
are the normalized associated Legendre functions defined
in~(\ref{association}),
$y_0$,~$y_1$, \dots, $y_{n-2}$,~$y_{n-1}$ are the positive zeros
of $\overline{P}^m_{m+2n+1}$ from~(\ref{zeros2}),
and $\sigma_0$,~$\sigma_1$, \dots, $\sigma_{n-2}$,~$\sigma_{n-1}$
are the corresponding quadrature weights from~(\ref{Christoffel_definition2}).
We will refer to the map via~(\ref{odds})
from $\nu_0$,~$\nu_1$, \dots, $\nu_{n-2}$,~$\nu_{n-1}$
to $\mu_0$,~$\mu_1$, \dots, $\mu_{n-2}$,~$\mu_{n-1}$
as the forward transform, and the map via~(\ref{odds})
from $\mu_0$,~$\mu_1$, \dots, $\mu_{n-2}$,~$\mu_{n-1}$
to $\nu_0$,~$\nu_1$, \dots, $\nu_{n-2}$,~$\nu_{n-1}$
as the inverse transform (as above, the inverse is also the transpose).

For the test vectors $\beta$ and $\nu$ whose entries appear
in~(\ref{evens}) and~(\ref{odds}), we used normalized vectors whose entries
were pseudorandom numbers drawn uniformly from $(-1,1)$,
normalized so that the sum of the squares of the entries is 1.

As described in Remark~\ref{adaptivity},
we compute for each block in the multilevel representation of $A$
an ID whose rank is as small as possible while still approximating
the block to nearly machine precision.

\begin{itemize}
\item[] The following list describes the headings of the tables:
\item $n$ is the size of the transform, the size of the vectors
      $\alpha$ and $\beta$ whose entries are given in~(\ref{evens}),
      and of the vectors $\mu$ and $\nu$ whose entries are given
      in~(\ref{odds}).
\item $m$ is the order of the normalized associated Legendre functions
      used in~(\ref{evens}) and~(\ref{odds}).
\item $k_{\rm max}$ is the maximum of the ranks of the IDs for the blocks
      in the multilevel representation.
\item $k_{\rm avg}$ is the average of the ranks of the IDs for the blocks
      in the multilevel representation.
\item $k_{\sigma}$ is the standard deviation of the ranks of the IDs
      for the blocks in the multilevel representation.
\item $t_{\rm dir}$ is the time in seconds required to apply an $n \times n$
      matrix to an $n \times 1$ vector using the standard procedure.
      We estimated the last two entries for $t_{\rm dir}$ by multiplying
      the third-to-last entry in each table by 4 and 16,
      since the large matrices required to generate those entries
      cannot fit in the available 2~GB of RAM. We indicate that these entries
      are estimates by enclosing them in parentheses.
\item $t_{\rm fwd}$ is the time in seconds required by the butterfly algorithm
      to compute the forward transform via~(\ref{evens}) or~(\ref{odds}),
      mapping from $\beta_0$,~$\beta_1$, \dots, $\beta_{n-2}$,~$\beta_{n-1}$
      to $\alpha_0$,~$\alpha_1$, \dots, $\alpha_{n-2}$,~$\alpha_{n-1}$,
      or from $\nu_0$,~$\nu_1$, \dots, $\nu_{n-2}$,~$\nu_{n-1}$
      to $\mu_0$,~$\mu_1$, \dots, $\mu_{n-2}$,~$\mu_{n-1}$.
\item $t_{\rm inv}$ is the time in seconds required by the butterfly algorithm
      to compute the inverse transform via~(\ref{evens}) or~(\ref{odds}),
      mapping from $\alpha_0$,~$\alpha_1$, \dots,
      $\alpha_{n-2}$,~$\alpha_{n-1}$
      to $\beta_0$,~$\beta_1$, \dots, $\beta_{n-2}$,~$\beta_{n-1}$,
      or from $\mu_0$,~$\mu_1$, \dots, $\mu_{n-2}$,~$\mu_{n-1}$
      to $\nu_0$,~$\nu_1$, \dots, $\nu_{n-2}$,~$\nu_{n-1}$.
\item $t_{\rm quad}$ is the time in seconds required in the precomputations
      to compute the quadrature nodes $x_0$,~$x_1$, \dots, $x_{n-2}$,~$x_{n-1}$
      or $y_0$,~$y_1$, \dots, $y_{n-2}$,~$y_{n-1}$,
      and weights $\rho_0$,~$\rho_1$, \dots, $\rho_{n-2}$,~$\rho_{n-1}$
      or $\sigma_0$,~$\sigma_1$, \dots, $\sigma_{n-2}$,~$\sigma_{n-1}$,
      from~(\ref{Christoffel_definition}) or~(\ref{Christoffel_definition2}),
      used in~(\ref{evens}) and~(\ref{odds}).
\item $t_{\rm comp}$ is the time in seconds required to construct
      the compressed multilevel representation used in the butterfly algorithm,
      after having already computed
      the quadrature nodes $x_0$,~$x_1$, \dots, $x_{n-2}$,~$x_{n-1}$
      or $y_0$,~$y_1$, \dots, $y_{n-2}$,~$y_{n-1}$,
      and weights $\rho_0$,~$\rho_1$, \dots, $\rho_{n-2}$,~$\rho_{n-1}$
      or $\sigma_0$,~$\sigma_1$, \dots, $\sigma_{n-2}$,~$\sigma_{n-1}$,
      from~(\ref{Christoffel_definition}) or~(\ref{Christoffel_definition2}),
      used in~(\ref{evens}) and~(\ref{odds}).
\item $m_{\rm max}$ is the maximum number of floating-point words
      of memory required to store entries of the transform matrix
      during any point in the precomputations
      (all other memory requirements are negligible in comparison).
\item $\epsilon_{\rm fwd}$ is the maximum difference between the entries
      in the result of the forward transform computed
      via the butterfly algorithm and those computed directly
      via the standard procedure for applying a matrix to a vector.
      (The result of the forward transform is the vector $\alpha$
      whose entries are given in~(\ref{evens})
      or the vector $\mu$ whose entries are given in~(\ref{odds}).)
\item $\epsilon_{\rm inv}$ is the maximum difference between the entries
      in a test vector and the entries in the result of applying
      to the test vector first the forward transform and
      then the inverse transform, both computed via the butterfly algorithm.
      (The result of the forward transform is the vector $\alpha$
      whose entries are given in~(\ref{evens})
      or the vector $\mu$ whose entries are given in~(\ref{odds}).
      The result of the inverse transform is the vector $\beta$
      whose entries are given in~(\ref{evens})
      or the vector $\nu$ whose entries are given in~(\ref{odds}).)
      Thus, $\epsilon_{\rm inv}$ measures the accuracy
      of the butterfly algorithm without reference to the standard procedure
      for applying a matrix to a vector (unlike $\epsilon_{\rm fwd}$).
\end{itemize}

For the first level of the multilevel representation
of the $n \times n$ matrix, we partitioned the matrix
into blocks each dimensioned $n \times 60$
(except for the rightmost block, since $n$ is not divisible by 60).
Every block on every level has about the same number of entries
(specifically, $60n$ entries).
We wrote all code in Fortran 77, compiling it using
the Lahey-Fujitsu Linux Express v6.2 compiler, with optimization flag
{\tt {-}{-}o2} enabled.
We ran all examples on one core of a 2.7~GHz Intel Core~2 Duo microprocessor
with 3~MB of L2 cache and 2~GB of RAM.
As described in Section~\ref{spharmonics},
we used extended-precision arithmetic during the portion of the preprocessing
requiring integration of an ordinary differential equation,
to compute quadrature nodes and weights
(this is not necessary to attain high accuracy, but does yield a couple
of extra digits of precision).
Otherwise, our code is compliant with the IEEE double-precision standard
(so that the mantissas of variables have approximately one bit of precision
less than 16 digits, yielding a relative precision of about .2E--15).

\begin{remark}
\label{empirical_rank}
Tables~\ref{evens_time}, \ref{small_time}, and~\ref{odds_time}
indicate that the linear transformations in~(\ref{evens}) and~(\ref{odds})
satisfy the rank property discussed in Subsection~\ref{butterfly},
with arbitrarily high precision, at least in some averaged sense.
Furthermore, it appears that the parameter $k$ discussed
in Subsection~\ref{butterfly} can be set to be independent
of the order $m$ and size $n$ of the transforms in~(\ref{evens})
and~(\ref{odds}), with the parameter $C$ discussed
in Subsection~\ref{butterfly} roughly proportional to $k$.
The acceleration provided by the butterfly algorithm thus
is sufficient for computing fast spherical harmonic transforms,
and is competitive with the approach taken in~\cite{tygert_sph}
(though much work remains in optimizing both approaches
in order to gauge their relative performance).
Unlike the approach taken in~\cite{tygert_sph},
the approach of the present paper does not require the use
of extended-precision arithmetic during the precomputations
in order to attain accuracy close to the machine precision,
even while accelerating spherical harmonic transforms about as well.
Moreover, the butterfly can be easier to implement.
\end{remark}

\begin{remark}
The values in Tables~\ref{evens_pre}, \ref{small_pre}, and~\ref{odds_pre}
vary with the size $n$ of the transforms in~(\ref{evens}) and~(\ref{odds})
as expected.
The values for $t_{\rm quad}$
are consistent with the expected values of a constant times $n$.
The values for $t_{\rm comp}$
are consistent with the expected values of a constant times $n^2$
(with the constant being proportional to $k_{\rm avg}$).
The values for $m_{\rm max}$
are consistent with the expected values of a constant times $n \log(n)$
(again with the constant being proportional to $k_{\rm avg}$);
these modest memory requirements make the preprocessing feasible
for large values of $n$ such as those in the tables.
\end{remark}

\begin{remark}
In the current technological environment,
neither the scheme of~\cite{tygert_sph} nor the approach of the present paper
is uniformly superior to the other.
For example, the theory from~\cite{tygert_sph} is rigorous
and essentially complete, while the theory of the present article ideally
should undergo further development, to prove that the rank properties discussed
in Remark~\ref{empirical_rank} are as strong as numerical experiments indicate,
yielding the desired acceleration.
In contrast, to attain accuracy close to the machine precision,
the approach of~\cite{tygert_sph} requires the use of extended-precision
arithmetic during its precomputations,
whereas the scheme of the present paper does not.
Implementing the procedure of the present article can be easier.
Finally, an anonymous referee kindly compared the running-times
of the implementations reported in~\cite{tygert_sph} and the present paper,
noticing that the newer computer system used in the present article is
about 2.7 times faster than the old system; the algorithms of~\cite{tygert_sph}
and of the present article are roughly equally efficient ---
certainly neither appears to be more than twice faster than the other.
However, both implementations are rather crude, and could undoubtedly benefit
from further optimization by experts on computer architectures;
also, we made no serious attempt to optimize the precomputations.
Furthermore, with the advent of multicore and distributed processors,
coming changes in computer architectures might affect
the two approaches differently, as they may parallelize and utilize cache
in different ways.
In the end, the use of one approach rather than the other may be
a matter of convenience, as the two methods are yielding similar performance.
\end{remark}

\section{Conclusions}
\label{conclusions}

This article provides an alternative means for performing the key
computational step required in~\cite{tygert_sph} for computing
fast spherical harmonic transforms.
Unlike the implementation described in~\cite{tygert_sph}
of divide-and-conquer spectral methods,
the butterfly scheme of the present paper does not require the use
of extended precision during the compression precomputations in order
to attain accuracy very close to the machine precision.
With the butterfly, the required amount of preprocessing is quite reasonable,
certainly not prohibitive.

Unfortunately, there seems to be little theoretical understanding
of why the butterfly procedure works so well for associated Legendre functions
(are the associated transforms nearly weighted averages
of Fourier integral operators?).
Complete proofs such as those in~\cite{oneil-woolfe-rokhlin}
and~\cite{tygert_sph} are not yet available for the scheme
of the present article.
By construction, the butterfly enables fast, accurate
applications of matrices to vectors when the precomputations succeed.
However, we have yet to prove that the precomputations
will compress the appropriate $n \times n$ matrix
enough to enable applications of the matrix to vectors
using only $\bigoh(n \log(n))$ floating-point operations (flops).
Nevertheless, the scheme has succeeded in all our numerical tests.
We hope to produce rigorous mathematical proofs that the precomputations
always compress the matrices as much as they did in our numerical experiments.

The precomputations for the algorithm of the present article
require $\bigoh(n^2)$ flops.
The precomputations for the algorithm of~\cite{tygert_sph} also
require $\bigoh(n^2)$ flops as implemented for the numerical examples
of that paper; however, the procedure of~\cite{gu-eisenstat95}
leads naturally to precomputations for the approach of~\cite{tygert_sph}
requiring only $\bigoh(n \log(n))$ flops
(though these ``more efficient'' precomputations do not become more efficient
in practice until $n$ is absurdly large, too large even to estimate reliably).
We do not expect to be able to accelerate the precomputations
for the algorithm of the present article without first producing
the rigorous mathematical proofs mentioned in the previous paragraph.
Even so, the current amount of preprocessing is not unreasonable,
as the numerical examples of Section~\ref{numerical} illustrate.

\section*{Acknowledgements}

We would like to thank V. Rokhlin for his advice, for his encouragement,
and for the use of his software libraries,
all of which have greatly enhanced this paper
and the associated computer codes.
We are also grateful to R. R. Coifman and Y. Shkolnisky.
We would like to thank the anonymous referees for their useful suggestions.

\begin{table}[p]
\caption{Ranks and running-times for even degrees}
\label{evens_time}
\vspace{1em}
\begin{tabular*}{\columnwidth}{@{\extracolsep{\fill}}rrcccccc}
  $n$ &   $m$ & $k_{\rm max}$ & $k_{\rm avg}$ & $k_{\sigma}$ & $t_{\rm dir}$ & $t_{\rm fwd}$ & $t_{\rm inv}$ \\\hline
 1250 &  1250 &           170 &          65.3 &         29.5 &       .25E--2 &       .18E--2 &       .15E--2 \\\hline
 2500 &  2500 &           168 &          67.0 &         32.6 &       .98E--2 &       .46E--2 &       .38E--2 \\\hline
 5000 &  5000 &           195 &          70.5 &         35.9 &       .39E--1 &       .12E--1 &       .96E--2 \\\hline
10000 & 10000 &           247 &          73.5 &         38.7 &       .15E--0 &       .28E--1 &       .23E--1 \\\hline
20000 & 20000 &           308 &          75.9 &         41.5 &      (.60E--0)&       .67E--1 &       .55E--1 \\\hline
40000 & 40000 &           379 &          78.0 &         43.8 &       (.24E+1)&       .16E--0 &       .13E--0 \\\hline
\end{tabular*}
\end{table}

\begin{table}[p]
\caption{Precomputation times, memory requirements, and accuracies
         for even degrees}
\label{evens_pre}
\vspace{1em}
\begin{tabular*}{\columnwidth}{@{\extracolsep{\fill}}rrccccc}
  $n$ &   $m$ & $t_{\rm quad}$ & $t_{\rm comp}$ & $m_{\rm max}$ & $\epsilon_{\rm fwd}$ & $\epsilon_{\rm inv}$ \\\hline
 1250 &  1250 &          .46E2 &          .96E0 &         .86E6 &             .62E--14 &             .19E--13 \\\hline
 2500 &  2500 &          .92E2 &          .41E1 &         .20E7 &             .37E--14 &             .25E--13 \\\hline
 5000 &  5000 &          .18E3 &          .17E2 &         .50E7 &             .59E--14 &             .43E--13 \\\hline
10000 & 10000 &          .37E3 &          .82E2 &         .14E8 &             .32E--14 &             .57E--13 \\\hline
20000 & 20000 &          .74E3 &          .39E3 &         .29E8 &             .30E--14 &             .88E--13 \\\hline
40000 & 40000 &          .15E4 &          .17E4 &         .64E8 &             .24E--14 &             .13E--12 \\\hline
\end{tabular*}
\end{table}

\begin{table}[p]
\caption{Ranks and running-times for even degrees}
\label{small_time}
\vspace{1em}
\begin{tabular*}{\columnwidth}{@{\extracolsep{\fill}}rrcccccc}
  $n$ &   $m$ & $k_{\rm max}$ & $k_{\rm avg}$ & $k_{\sigma}$ & $t_{\rm dir}$ & $t_{\rm fwd}$ & $t_{\rm inv}$ \\\hline
 1250 &     0 &           110 &          67.0 &         23.3 &       .25E--2 &       .18E--2 &       .15E--2 \\\hline
 2500 &     0 &           110 &          70.0 &         25.0 &       .98E--2 &       .48E--2 &       .40E--2 \\\hline
 5000 &     0 &           111 &          73.9 &         26.1 &       .39E--1 &       .12E--1 &       .10E--1 \\\hline
10000 &     0 &           111 &          77.3 &         26.8 &       .15E--0 &       .29E--1 &       .24E--1 \\\hline
20000 &     0 &           112 &          80.2 &         27.1 &      (.60E--0)&       .68E--1 &       .57E--1 \\\hline
40000 &     0 &           169 &          82.7 &         27.4 &       (.24E+1)&       .16E--0 &       .13E--0 \\\hline
\end{tabular*}
\end{table}

\begin{table}[p]
\caption{Precomputation times, memory requirements, and accuracies
         for even degrees}
\label{small_pre}
\vspace{1em}
\begin{tabular*}{\columnwidth}{@{\extracolsep{\fill}}rrccccc}
  $n$ &   $m$ & $t_{\rm quad}$ & $t_{\rm comp}$ & $m_{\rm max}$ & $\epsilon_{\rm fwd}$ & $\epsilon_{\rm inv}$ \\\hline
 1250 &     0 &          .44E2 &          .96E0 &         .86E6 &             .49E--14 &             .12E--12 \\\hline
 2500 &     0 &          .88E2 &          .41E1 &         .20E7 &             .35E--14 &             .14E--12 \\\hline
 5000 &     0 &          .18E3 &          .18E2 &         .51E7 &             .23E--14 &             .35E--12 \\\hline
10000 &     0 &          .36E3 &          .82E2 &         .14E8 &             .18E--14 &             .63E--12 \\\hline
20000 &     0 &          .74E3 &          .40E3 &         .29E8 &             .20E--14 &             .22E--11 \\\hline
40000 &     0 &          .14E4 &          .18E4 &         .66E8 &             .16E--14 &             .37E--11 \\\hline
\end{tabular*}
\end{table}

\newpage

\begin{table}[h]
\caption{Ranks and running-times for odd degrees}
\label{odds_time}
\vspace{1em}
\begin{tabular*}{\columnwidth}{@{\extracolsep{\fill}}rrcccccc}
  $n$ &   $m$ & $k_{\rm max}$ & $k_{\rm avg}$ & $k_{\sigma}$ & $t_{\rm dir}$ & $t_{\rm fwd}$ & $t_{\rm inv}$ \\\hline
 1250 &  1250 &           170 &          65.3 &         29.5 &       .25E--2 &       .17E--2 &       .15E--2 \\\hline
 2500 &  2500 &           169 &          67.0 &         32.6 &       .98E--2 &       .48E--2 &       .39E--2 \\\hline
 5000 &  5000 &           196 &          70.5 &         35.9 &       .39E--1 &       .12E--1 &       .97E--2 \\\hline
10000 & 10000 &           247 &          73.5 &         38.7 &       .15E--0 &       .28E--1 &       .24E--1 \\\hline
20000 & 20000 &           308 &          75.9 &         41.4 &      (.60E--0)&       .68E--1 &       .56E--1 \\\hline
40000 & 40000 &           379 &          78.0 &         43.8 &       (.24E+1)&       .16E--0 &       .13E--0 \\\hline
\end{tabular*}
\end{table}

\begin{table}[h]
\caption{Precomputation times, memory requirements, and accuracies
         for odd degrees}
\label{odds_pre}
\vspace{1em}
\begin{tabular*}{\columnwidth}{@{\extracolsep{\fill}}rrccccc}
  $n$ &   $m$ & $t_{\rm quad}$ & $t_{\rm comp}$ & $m_{\rm max}$ & $\epsilon_{\rm fwd}$ & $\epsilon_{\rm inv}$ \\\hline
 1250 &  1250 &          .46E2 &          .94E0 &         .86E6 &             .41E--14 &             .19E--13 \\\hline
 2500 &  2500 &          .92E2 &          .40E1 &         .20E7 &             .41E--14 &             .29E--13 \\\hline
 5000 &  5000 &          .18E3 &          .17E2 &         .50E7 &             .40E--14 &             .51E--13 \\\hline
10000 & 10000 &          .37E3 &          .80E2 &         .14E8 &             .31E--14 &             .62E--13 \\\hline
20000 & 20000 &          .74E3 &          .39E3 &         .29E8 &             .34E--14 &             .10E--12 \\\hline
40000 & 40000 &          .15E4 &          .18E4 &         .64E8 &             .25E--14 &             .14E--12 \\\hline
\end{tabular*}
\end{table}

\vspace{.333in}





\section*{References}
\bibliographystyle{elsarticle-num_sorted}
\bibliography{but.bib}

\begin{thebibliography}{10}
\expandafter\ifx\csname url\endcsname\relax
  \def\url#1{\texttt{#1}}\fi
\expandafter\ifx\csname urlprefix\endcsname\relax\def\urlprefix{URL }\fi
\expandafter\ifx\csname href\endcsname\relax
  \def\href#1#2{#2} \def\path#1{#1}\fi

\bibitem{abramowitz-stegun}
M.~Abramowitz, I.~A. Stegun (Eds.), Handbook of Mathematical Functions, Dover
  Publications, New York, 1972.

\bibitem{adams-swarztrauber}
J.~C. Adams, P.~N. Swarztrauber, {SPHEREPACK} 3.0: {A} model development
  facility, Mon. Wea. Rev. 127 (1999) 1872--1878.

\bibitem{aho-hopcroft-ullman}
A.~Aho, J.~Hopcroft, J.~Ullman, Data Structures and Algorithms, Addison-Wesley,
  1987.

\bibitem{candes-demanet-ying}
E.~Cand\`es, L.~Demanet, L.~Ying, A fast butterfly algorithm for the
  computation of {F}ourier integral operators, Multiscale Model. Simul. 7~(4)
  (2009) 1727--1750.

\bibitem{cheng-gimbutas-martinsson-rokhlin}
H.~Cheng, Z.~Gimbutas, P.-G. Martinsson, V.~Rokhlin, On the compression of low
  rank matrices, SIAM J. Sci. Comput. 26~(4) (2005) 1389--1404.

\bibitem{goreinov-tyrtyshnikov}
S.~A. Goreinov, E.~E. Tyrtyshnikov, The maximal-volume concept in approximation
  by low-rank matrices, in: V.~Olshevsky (Ed.), Structured Matrices in
  Mathematics, Computer Science, and Engineering I: Proceedings of an
  AMS-IMS-SIAM Joint Summer Research Conference, University of Colorado,
  Boulder, June 27--July 1, 1999, Vol. 280 of Contemporary Mathematics, AMS
  Publications, Providence, RI, 2001, pp. 47--51.

\bibitem{gu-eisenstat95}
M.~Gu, S.~C. Eisenstat, A divide-and-conquer algorithm for the symmetric
  tridiagonal eigenproblem, SIAM J. Matrix Anal. Appl. 16 (1995) 172--191.

\bibitem{gu-eisenstat96}
M.~Gu, S.~C. Eisenstat, Efficient algorithms for computing a strong
  rank-revealing {QR} factorization, SIAM J. Sci. Comput. 17~(4) (1996)
  848--869.

\bibitem{martinsson-rokhlin-tygert1}
P.-G. Martinsson, V.~Rokhlin, M.~Tygert, On interpolation and integration in
  finite-dimensional spaces of bounded functions, Comm. Appl. Math. Comput.
  Sci. 1 (2006) 133--142.

\bibitem{michielssen-boag}
E.~Michielssen, A.~Boag, A multilevel matrix decomposition algorithm for
  analyzing scattering from large structures, IEEE Trans. Antennas and
  Propagation 44~(8) (1996) 1086--1093.

\bibitem{oneil-woolfe-rokhlin}
M.~O'Neil, F.~Woolfe, V.~Rokhlin, An algorithm for the rapid evaluation of
  special function transforms, Appl. Comput. Harmon. Anal. 28~(2) (2010)
  203--226.

\bibitem{reuter-ratner-seideman}
M.~G. Reuter, M.~A. Ratner, T.~Seideman, A fast method for solving both the
  time-dependent {S}chr\"odinger equation in angular coordinates and its
  associated ``$m$-mixing'' problem, J. Chem. Phys. 131 (2009) 094108--1 ---
  094108--6.

\bibitem{swarztrauber-spotz}
P.~N. Swarztrauber, W.~F. Spotz, Generalized discrete spherical harmonic
  transforms, J. Comput. Phys. 159~(2) (2000) 213--230.

\bibitem{szego}
G.~Szeg\"o, Orthogonal Polynomials, 11th Edition, Vol.~23 of Colloquium
  Publications, American Mathematical Society, Providence, RI, 2003.

\bibitem{tygert_sph}
M.~Tygert, Fast algorithms for spherical harmonic expansions, {II}, J. Comput.
  Phys. 227~(8) (2008) 4260--4279.

\bibitem{tyrtyshnikov}
E.~E. Tyrtyshnikov, Incomplete cross approximation in the mosaic-skeleton
  method, Computing 64~(4) (2000) 367--380.

\bibitem{ying}
L.~Ying, Sparse {F}ourier transform via butterfly algorithm, SIAM J. Sci.
  Comput. 31~(3) (2009) 1678--1694.

\end{thebibliography}







\end{document}